\documentclass[11pt,reqno]{amsart}
\usepackage[colorlinks=true, linkcolor=blue, urlcolor=blue]{hyperref}
\usepackage{graphicx}
\usepackage{url}
\usepackage{amsmath,amsthm,amssymb,amsbsy,amstext,amsfonts,amscd}
\usepackage{mathrsfs, color}

\newtheorem{thm}{Theorem}[section]

\newtheorem{lem}[thm]{Lemma}

\theoremstyle{definition}

\theoremstyle{remark}
\newtheorem{rem}{Remark}
\newtheorem{exam}{Example}

\numberwithin{equation}{section}

\begin{document}

\title[logarithmic integrals]{On  three classes of logarithmic integrals}
\author{Necdet Bat{\i}r}
\address{Department of  Mathematics, Nev{\c{s}}ehir Hac{\i} Bekta\c{ş} Veli University, Nev{\c{s}}ehir, 50300 Turkey}
\email{nbatir@hotmail.com}

\author{Nandan Sai Dasireddy}
\address{House No. 3-11-363, Pavani Plaza, Flat No. 401, Road No. 2, Shivaganga Colony Siris Road,\\
L. B. Nagar, Hyderabad, Telangana, India}
\email{dasireddy.1818@gmail.com}


\begin{abstract}
In this paper, we evaluate  the following families of definite integrals in closed form and we show that they are expressible only in terms of the dilogarithm function and the inverse tangent integral, and elementary functions.
\begin{equation*}
\int_{0}^{1}\frac{\log\big(x^m+1\big)}{x+1}\thinspace{\rm d}x \quad \mbox{and}\quad \int_{0}^{1}\frac{\log\big(x^m+1\big)}{x^2+1}\thinspace{\rm d}x,
\end{equation*}
where $m$ is a positive odd integer.  When $m$ is a positive even integer, these integrals have been evaluated previously by Sofo and Bat{\i}r, and the case where $m$ is an odd integer has been left as open problems. The integrals of the first kind arise in Zagier's work on the Kronecker limit formula. In addition, we demonstrate that a functional equation satisfied by the Herglotz-Zagier-Novikov function is a very specific case of  of a more general formula, and give  numerous illustrative examples.
\end{abstract}
\maketitle

\section{Introduction }\label{sec1}
Zagier's article \cite{Zagier 2} on the Kronecker limit formula has been a source of inspiration for many researchers working in number theory. One of them, Novikov \cite{Novikov} obtained a Kronecker limit formula involving  the following logarithmic integral
\begin{equation}\label{e01}
\mathscr{F}(x;u,v) =\int_{0}^{1}\frac{\log\big(1-ut^x\big)}{v^{-1}-t}{\rm d}t,
\end{equation}
where $\Re (x)>0$, $u\in\mathbb{C}\backslash{(1,\infty)}$, $u\neq1$ and $v\in\mathbb{C}\backslash[1,\infty)$, $v\neq1$. Choie and Kumar \cite{Choie} named this function the Herglotz-Zagier-Novikov function and discovered many of its properties. As an example, they proved that
\begin{equation}\label{e02}
\mathscr{F}(x;u,v) +\mathscr{F}\left(\frac{1}{x};v,u\right)=-\log(1-u)\log(1-v);
\end{equation}
see \cite[Theorem 2.3]{Choie}.  The function given in \eqref{e01} is a generalization of some functions previously studied by Herglotz \cite{Herglotz}, Radchenko and Zagier \cite{Zagier}, and Muzaffar and Williams \cite{Muzaf}. These functions are as follows:
\begin{equation*}
F(x)=\sum_{n=1}^{\infty}\frac{\psi(nx)-\log(nx)}{n},
\end{equation*}
where $\psi(x)=\frac{\Gamma^\prime(x)}{\Gamma(x)}$ is the digamma function,
\begin{equation*}
J(x)=\int_{0}^{1}\frac{\log\big(1+t^x\big)}{1+t}{\rm d}t
\end{equation*}
and
\begin{equation*}
T(x)=\int_{0}^{1}\frac{\arctan\big(t^x\big)}{1+t^2}{\rm d}t
\end{equation*}
for $\Re(x)>0$. Very recently Radchenko and Zagier \cite{Zagier} studied $F(x)$ and $J(x)$, and they found the following connection between them:
\begin{equation*}
  J(x)=F(2x)-2F(x)+F\left(\frac{x}{2}\right)+\frac{\pi^2}{12x}.
\end{equation*}
As special cases of $\mathscr{F}(x;u,v)$, by \eqref{e02} $J(x)$ and $T(x)$ satisfy the following simple functional relations:

\begin{equation}\label{e03}
J(x)+J\left(\frac{1}{x}\right)=\log^22
\end{equation}
and
\begin{equation*}
T(x)+T\left(\frac{1}{x}\right)=\frac{\pi^2}{16}.
\end{equation*}
In \cite{Zagier} and \cite{Choie} the authors evaluated $J(x)$ for many particular values of $x$. For example, in \cite{Zagier} the authors showed that
\begin{equation*}
J\big(4+\sqrt{17}\thinspace\big)=-\frac{\pi^2}{6}+\frac{1}{2}\log^22+\log(2)\log\big(4+\sqrt{17}\big)
\end{equation*}
and
\begin{equation*}
J\left(\frac{2}{5}\right)=\frac{11\pi^2}{240}+\frac{3}{4}\log^22-2\log\bigg(\frac{\sqrt{5}+1}{2}\bigg).
\end{equation*}
In \cite[p. 411]{So-Bat} Sofo and Bat{\i}r established the following formula
\begin{align}\label{e06}
\int_{0}^{1}\frac{\log\big(1+u^{2b}\big)}{1+u}{\rm d}u&=\frac{1}{2}\sum_{k=0}^{2b-1}\log^2\bigg(2\sin\left(\frac{(2k+1)\pi}{4b}\right)\bigg)+\frac{1-2b^2}{8b}\zeta(2)\notag\\
&+\log^22-\frac{1}{2}\sum_{k=0}^{b-1}\log^2\bigg(2\sin\left(\frac{(2k+1)\pi}{2b}\right)\bigg), \quad b \in \mathbb{N}
\end{align}
and in \cite[Example 17]{Bat-2} Bat{\i}r proved the following integral formula
\begin{equation*}
\int_{0}^{1}\frac{\log\big(x^3+1\big)}{x+1}{\rm d}x=\frac{1}{2}\operatorname{Li}_2\left(-\frac{1}{3}\right)+\frac{1}{2}\log^23-\frac{1}{6}\zeta(2)
\end{equation*}
and in \cite[Remark 7]{Bat-2} he  proposed the following problem, as an interesting area of source to study: Whether there exists an explicit computational formula for the logarithmic integral
\begin{equation*}
\int_{0}^{1}\frac{\log\big(x^m+1\big)}{x+1}{\rm d}x,
\end{equation*}
where $m$ is an odd integer  with $m\geq 5$. We note that a formula equivalent to \eqref{e06} has been achieved in \cite{Choie}. In 2023 Choie and Kumar \cite[Corollary 3.5]{Choie} provided the following  general formula for $J(n)$, which is valid for all $n\in \mathbb{N}$:
\begin{equation*}
J(n)=\frac{n}{2}\log^2 2-\frac{\pi^2(n^2-1)}{12n}+\sum_{j=1}^{n}\operatorname{Li}_2\bigg(\frac{1}{2}\left(1+e^{\frac{\pi i(2j+1)}{n}}\right)\bigg).
\end{equation*}
This formula looks very nice, but it does not give explicit results even for small values of positive odd integers $n$. For example, for $n=3$ we have
\begin{equation*}
J(3)=\frac{3}{2}\log^2 2-\frac{2\pi^2}{9}+2\,\Re\bigg\{\operatorname{Li}_2\bigg(\frac{3+i\sqrt{3}}{4}\thinspace\bigg)\bigg\}
\end{equation*}
and it is not easy to derive formula \eqref{e3.4}  from this identity. Our first aim in this work is to evaluate $J(m)$ for all odd integers $m$  in closed form in terms of real arguments of the dilogarithm function $\operatorname{Li}_2(x)$.

In \cite[Theorem 2.5]{Bat-2} Bat{\i}r has evaluated the following family of integrals in closed form for even positive integers $m$:
\begin{equation*}
\int_{0}^{1}\frac{\log\big(x^m+1\big)}{x^2+1}{\rm d}x.
\end{equation*}
He also showed that
\begin{equation}\label{e08.1}
\int_{0}^{1}\frac{\log\big(x^3+1\big)}{x^2+1}{\rm d}x=\frac{\pi}{8}\log 2-\frac{5}{3}\operatorname{G}+\frac{\pi}{3}\log\big(2+\sqrt{3}\big)
\end{equation}
and left the evaluation of these integrals for odd integers $m$ with $m\geq 5$ as an open problem. Our second aim is to provide a solution to this problem and to evaluate these integrals in closed form for all odd integers $m\geq 3$.

Our third and last aim is to provide a generalization of the Herglotz-Zagier-Novikov function given by \eqref{e01} and to offer some
applications.  In this work, an empty sum is conventionally interpreted as zero wherever it appears.

We need the following three lemmas.
\begin{lem}\label{lem1}Let $m$ be an odd integer with $m\geq3$. Then we have
\begin{align*}
\log\bigg[\bigg(\frac{1-u}{1+u}\bigg)^m+1\bigg]=&\frac{m+1}{2}\log 2-m\log(1+u)+\sum_{k=0}^{\frac{m-3}{2}}\log(1+\varphi_k)\\
+&\sum_{k=0}^{\frac{m-3}{2}}\log\bigg[u^2+\frac{1-\varphi_k}{1+\varphi_k}\bigg],
\end{align*}
where $\varphi_k=\cos\left(\frac{(2k+1)\pi}{m}\right)$.
\end{lem}
\begin{proof}
By \cite[p.\,5]{Bat-2} we have
\begin{align*}
z^m+1&=(z+1)\prod_{k=0}^{\frac{m-3}{2}}\bigg(z-e^{\frac{ i(2k+1)\pi}{m}}\bigg)\bigg(z-e^{-\frac{ i(2k+1)\pi}{m}}\bigg)\\
&=(z+1)\prod_{k=0}^{\frac{m-3}{2}}\big[z^2-2z\varphi_k+1\big].
\end{align*}
Setting $z=\frac{1-u}{1+u}$, and then taking the logarithm of both sides, we, after an easy computation, get the desired result.
\end{proof}


\begin{lem}\label{lem2}Let $q$ be any nonnegative real number. Then we have
\begin{align*}
\int_{0}^{1}\frac{\log\big(u^2+q\big)}{1+u}{\rm d}u=\log 2\log(1+q)-\arctan^2\left(\frac{1}{\sqrt{q}}\right)+\frac{1}{2}\operatorname{Li}_2\left(\frac{1}{q+1}\right).
\end{align*}
\end{lem}
\begin{proof}
By the substitution $ u=t\sqrt{q}$, we have
\begin{equation*}
\int_{0}^{1}\frac{\log\big(u^2+q\big)}{u+1}{\rm d}u=\sqrt{q}\log q\int_{0}^{1/\sqrt{q}}\frac{{\rm d}t}{1+\sqrt{q}t}+\sqrt{q}\int_{0}^{1/\sqrt{q}}\frac{\log\big(t^2+1\big)}{1+t\sqrt{q}}{\rm d}t.
\end{equation*}
Evaluating the first integral and performing integration by parts on the second yields
\begin{equation*}
\int_{0}^{1}\frac{\log\big(u^2+q\big)}{u+1}{\rm d}u=\log 2\log(q+1)-2\int_{0}^{1/\sqrt{q}}\frac{t\log\big(1+t\sqrt{q}\big)}{1+t^2}{\rm d}t.
\end{equation*}

Setting $t=\frac{v}{\sqrt{q}}$, we obtain
\begin{equation*}
\int_{0}^{1}\frac{\log\big(u^2+q\big)}{u+1}{\rm d}u=\log 2\log(q+1)-\frac{2}{q}\int_{0}^{1}\frac{v\log(1+v)}{1+v^2/q}{\rm d}v.
\end{equation*}
The following result appears in  \cite[Eq. (1.7), p. 3]{Valean2}:
\begin{equation*}
\int_{0}^{1}\frac{x\log(1+x)}{1+ax^2}{\rm d}x=\frac{1}{2a}\arctan^2\big(\sqrt{a}\,\big)-\frac{1}{4a}\operatorname{Li}_2\left(\frac{a}{a+1}\right)\, (a>0).
\end{equation*}
From this result with $a=\frac{1}{q}$ it follows that
\begin{equation*}
\int_{0}^{1}\frac{\log\big(u^2+q\big)}{1+u}{\rm d}u=\log 2\log(1+q)-\arctan^2\left(\frac{1}{\sqrt{q}}\right)+\frac{1}{2}\operatorname{Li}_2\left(\frac{1}{q+1}\right).
\end{equation*}

\end{proof}


\begin{lem}\label{lem3}Let $q$ be a positive real number. Then we have
\begin{align}\label{e10}
\int_{0}^{1}\frac{\log \big(x^2+q\big)}{1+x^2}{\rm d}x=\frac{\pi}{2}\log\big(1+\sqrt{q}\big)+\operatorname{Ti}_2\left(\frac{\sqrt{q}-1}{\sqrt{q}+1}\right)-\operatorname{G},
\end{align}
where $\operatorname{G}$ is the Catalan's constant defined by
\begin{equation*}
\operatorname{G}=\sum_{k=0}^{\infty}\frac{(-1)^k}{(2k+1)^2}=0.915965...,
\end{equation*}
and $\operatorname{Ti}_2$ is the inverse tangent integral defined by
\begin{equation}\label{e10.1}
\operatorname{Ti}_2(x)=\int_{0}^{x}\frac{\arctan u}{u}{\rm d}u.
\end{equation}
 \end{lem}
\begin{proof}
We define
\begin{equation*}
g(q)=\int_{0}^{1}\frac{\log\big(x^2+q\big)}{x^2+1}{\rm d}x.
\end{equation*}
Differentiation with respect to $q$ gives
\begin{align*}
g^\prime(q)&=\int_{0}^{1}\frac{{\rm d}x}{\big(x^2+q\big)\big(x^2+1\big)}\notag\\
&=\frac{1}{q-1}\int_{0}^1\frac{{\rm d}x}{x^2+1}-\frac{1}{q-1}\int_{0}^1\frac{{\rm d}x}{x^2+q}\notag\\
&=\frac{\pi}{4}\frac{1}{q-1}-\frac{1}{q-1}\int_{0}^1\frac{{\rm d}x}{x^2+q}.
\end{align*}
It is very easy to see
\begin{align*}
\int_{0}^1\frac{{\rm d}x}{x^2+q}=\frac{\arctan\left(\frac{1}{\sqrt{q}}\right)}{\sqrt{q}}.
\end{align*}
Thus, we have
\begin{align*}
g^\prime(q)&=\frac{\pi}{4}\frac{1}{q-1}-\frac{\arctan\left(\frac{1}{\sqrt{q}}\right)}{(q-1)\sqrt{q}}.
\end{align*}
Integration gives
\begin{align*}
g(q)&=\frac{\pi}{4}\log|q-1|-\int\frac{\arctan\left(\frac{1}{\sqrt{q}}\right)}{(q-1)\sqrt{q}}{\rm d}q+C.
\end{align*}
Making the change of variable $u=\frac{1}{\sqrt{q}}$  we see that
\begin{align*}
g(q)&=\frac{\pi}{4}\log|q-1|+2\int\frac{\arctan(u)}{1-u^2}{\rm d}u+C\\
&=\frac{\pi}{4}\log|q-1|+2\int_{0}^{u}\frac{\arctan(y)}{1-y^2}{\rm d}y+C\\
&=\frac{\pi}{4}\log|q-1|-\int_{0}^{u}\left(-\frac{1}{1-y}-\frac{1}{1+y}\right)\arctan(y)\thinspace{\rm d}y+C.
\end{align*}
Applying integration by parts we obtain
\begin{align*}
g(q)&=\frac{\pi}{4}\log|q-1|-\log\left|\frac{1-u}{1+u}\right|\arctan(u)+\int_{0}^{u}\frac{1}{y^2+1}\log\bigg|\frac{1-y}{1+y}\bigg|{\rm d}y+C.
\end{align*}
Letting $q\to 1$ yields
\begin{align}\label{e11}
\int_{0}^{1}\frac{\log\big(x^2+1\big)}{x^2+1}{\rm d}x=&\lim\limits_{q\to 1}\bigg[\frac{\pi}{4}\log|q-1|-\log\frac{|1-u|}{1+u}\arctan (u)\bigg]\notag\\
&+\int_{0}^{1}\frac{1}{y^2+1}\log\bigg|\frac{1-y}{1+y}\bigg|{\rm d}y+C \quad \big(u=\frac{1}{\sqrt{q}}\big).
\end{align}
By \cite[Example 1]{Bat-2} (see also \cite[Lemma 1.5]{Nantan}) we have
\begin{align}\label{e12}
\int_{0}^{1}\frac{\log\big(x^2+1\big)}{x^2+1}{\rm d}x=\frac{\pi}{2}\log 2-\operatorname{G}.
\end{align}
We can also easily evaluate that
\begin{align}\label{e13}
\lim\limits_{q\to 1}\bigg[\frac{\pi}{4}\log|q-1|-\log\left|\frac{1-u}{1+u}\right|\arctan(u)\bigg]=\frac{\pi}{4}\log 2
\end{align}
and
\begin{align}\label{e14}
\int_{0}^{1}\frac{1}{y^2+1}\log\left|\frac{1-u}{1+u}\right|{\rm d}y=-\operatorname{G}.
\end{align}
Combining \eqref{e11}, \eqref{e12}, \eqref{e13} and \eqref{e14} we conclude that $C=0$. Thus we have

\begin{align}\label{e15}
g(q)&=\frac{\pi}{4}\log|q-1|-\log\left|\frac{1-u}{1+u}\right|\arctan(u)+\int_{0}^{u}\frac{1}{1+y^2}\log\bigg|\frac{1-y}{1+y}\bigg|{\rm d}y.
\end{align}
We recall the following two integral formulas from \cite[Eqs. (1.22) and (1.23)]{Valean2}:
\begin{align*}
\int_{0}^{a}\frac{\log (1-y)}{1+y^2}{\rm d}y=&\arctan a\log(1-a)+\frac{\pi}{4}\mbox{arctanh}\,a -\frac{1}{2}\operatorname{G}-\frac{1}{2}\operatorname{Ti}_2(a)\\
&+\frac{1}{2}\operatorname{Ti}_2\left(\frac{1-a}{1+a}\right)+\frac{1}{4}\operatorname{Ti}_2\left(\frac{2a}{1-a^2}\right)\quad |a|<1
\end{align*}
and
\begin{align*}
\int_{0}^{a}\frac{\log (1+y)}{1+y^2}{\rm d}y=&\arctan(a)\log(1+a)-\frac{\pi}{4}\mbox{arctanh}(a)  +\frac{1}{2}\operatorname{G}-\frac{1}{2}\operatorname{Ti}_2(a)\\
&-\frac{1}{2}\operatorname{Ti}_2\left(\frac{1-a}{1+a}\right)+\frac{1}{4}\operatorname{Ti}_2\left(\frac{2a}{1-a^2}\right)\quad |a|<1.
\end{align*}
Substituting the values of these two integrals above in \eqref{e15}, after replacing $a$ by $1/\sqrt{q}$, we find that
\begin{align*}
\int_{0}^{1}\frac{\log \big(x^2+q\big)}{1+x^2}{\rm d}x=\frac{\pi}{2}\log\big(1+\sqrt{q}\big)+\operatorname{Ti}_2\left(\frac{\sqrt{q}-1}{\sqrt{q}+1}\right)-\operatorname{G}.
\end{align*}
\end{proof}

\begin{rem}We would like to express our sincere thanks to the Editor for bringing to our attention the paper by \cite{Rogers}, in which it is recorded that $\operatorname{Ti}_2(x)=\Im\left(\operatorname{Li}_2(ix)\right)$ for $|x|\leq 1$.  By using the following functional relation
$$
\operatorname{Ti}_2(x)-\operatorname{Ti}_2(1/x)=\frac{\pi}{2}\log x
$$
for $\Re(x)>0$, the value of the integral in Lemma 1.2 can be expressed in terms of comples arguments of  the dilogarithm function but we prefer leaving it in its present form. Although $\operatorname{Ti}_2$ is not well known as the others, it has been studied by Ramanujan in \cite{Ramanujan} and many of its properties has been discovered. It is worth noting that Entry 4.295.22 in \cite[p. 561]{Grad} provides a related generalization, from which the following integral formula can be obtained:
\begin{equation*}
\int_{0}^{\infty}\frac{\log\big(x^2+q\big)}{x^2+1}{\rm d}x=\pi\log\big(1+\sqrt{q}\big).
\end{equation*}
\end{rem}


\section{Main results}\label{sec2}
We begin this section by presenting a solution to the problem proposed in \cite[Remark 7]{Bat-2}.
\begin{thm}\label{thm2.1} Let $m$ be any positive odd integer. Then
\begin{align}\label{e1}
J(m)=\frac{m}{2}\log^22-\frac{(m^2-1)\pi^2}{24m}+\frac{1}{2}\sum_{k=0}^{\frac{m-3}{2}}\operatorname{Li}_2\bigg(\cos^2\left(\frac{(2k+1)\pi}{2m}\right)\bigg).
\end{align}
Also, in view of \eqref{e03}, we have
\begin{align*}
J\left(\frac{1}{m}\right)=\frac{2-m}{2}\log^22+\frac{(m^2-1)\pi^2}{24m}-\frac{1}{2}\sum_{k=0}^{\frac{m-3}{2}}\operatorname{Li}_2\bigg(\cos^2\left(\frac{(2k+1)\pi}{2m}\right)\bigg),
\end{align*}
where $\varphi_k=\cos\left(\frac{(2k+1)\pi}{m}\right)$, and $\operatorname{Li}_2$ is the dilogarithm function defined by
\begin{equation*}
\operatorname{Li}_2(z)=\sum_{k=1}^{\infty}\frac{z^k}{k^2}\quad (|z|\leq 1).
\end{equation*}
\end{thm}
\begin{proof}
We make the substitution $x=\frac{1-u}{1+u}$ and ${\rm d}x=\frac{-2}{(1+u)^2}{\rm d}u$. It follows that
\begin{align*}
J(m)=\int_{0}^{1}\frac{1}{u+1}\log\left[\left(\frac{1-u}{1+u}\right)^m+1\right]{\rm d}u.
\end{align*}
By Lemma \ref{lem1} we get
\begin{align*}
J(m)=&\frac{1}{2}\log^22+\log2\sum_{k=0}^{\frac{m-3}{2}}\log(1+\varphi_k)+\sum_{k=0}^{\frac{m-3}{2}}\int_{0}^{1}\frac{1}{u+1}\log\bigg[u^2+\frac{1-\varphi_k}{1+\varphi_k}\bigg]{\rm d}u.
\end{align*}
Applying Lemma \ref{lem2} it follows that
\begin{align*}
J(m)=&\frac{m}{2}\log^22+\log2\sum_{k=0}^{\frac{m-3}{2}}\log(1+\varphi_k)+\sum_{k=0}^{\frac{m-3}{2}}\bigg[\log2\log\bigg(\frac{2}{1+\varphi_k}\bigg)\\
&-\arctan^2\bigg(\sqrt{\frac{1+\varphi_k}{1-\varphi_k}}\thinspace\bigg)+\frac{1}{2}\operatorname{Li}_2\bigg(\frac{1+\varphi_k}{2}\bigg)\bigg].
\end{align*}
Simplifying this identity, we get
\begin{align}\label{e02.1}
J(m)=\frac{m}{2}\log^22-\sum_{k=0}^{\frac{m-3}{2}}\arctan^2\bigg(\sqrt{\frac{1+\varphi_k}{1-\varphi_k}}\thinspace\bigg)+\frac{1}{2}\sum_{k=0}^{\frac{m-3}{2}}\operatorname{Li}_2\bigg(\frac{1+\varphi_k}{2}\bigg).
\end{align}
Since $1+\cos x=2\cos^2\frac{x}{2}$ and as can be easily shown that
\begin{equation*}
\arctan\left(\sqrt{\frac{1+\varphi_k}{1-\varphi_k}}\thinspace\right)=\frac{\pi}{4}+\frac{1}{2}\arcsin(\varphi_k)
\end{equation*}
and
\begin{equation*}
\arcsin(\varphi_k)=\frac{\pi}{2}-\frac{(2k+1)\pi}{m},
\end{equation*}
we arrive at
\begin{equation*}
\arctan\left(\sqrt{\frac{1+\varphi_k}{1-\varphi_k}}\thinspace\right)=\frac{(2k+1)\pi}{2m}.
\end{equation*}
Replacing this in \eqref{e02.1}, we get the result.
\end{proof}


\begin{thm}\label{thm2.2}Let $m$ be any positive odd integer. Then
\begin{align*}
&\int_{0}^{1}\frac{\log\big(x^m+1\big)}{x^2+1}{\rm d}x=\frac{m\pi}{8}\log2-\frac{(m-1)\operatorname{G}}{2}\notag\\
&+\frac{\pi}{4}\sum_{k=0}^{\frac{m-3}{2}}\log\left(1+\sin\frac{(2k+1)\pi}{m}\right)
+\sum_{k=0}^{\frac{m-3}{2}}\operatorname{Ti}_2\bigg(\frac{1}{\varphi_k}\left(\sqrt{1-\varphi_k^2}-1\right)\bigg),
\end{align*}
where $\operatorname{Ti}_2$ is the inverse tangent  integral  defined in \eqref{e10.1}.
\end{thm}
\begin{proof}
The proof is  based on the same arguments used in the proof of Theorem \ref{thm2.1}.  We define
\begin{equation*}
I(m):=\int_{0}^{1}\frac{\log\big(x^m+1\big)}{x^2+1}{\rm d}x.
\end{equation*}
Making the change of variable $x=\frac{1-u}{1+u}$, so that,  ${\rm d}x=-\frac{2}{(1+u)^2}{\rm d}u$, we get
\begin{align*}
I(m)=\int_{0}^{1}\log\bigg(1+\left(\frac{1-u}{1+u}\right)^m\bigg)\frac{{\rm d}u}{u^2+1}.
\end{align*}
Applying Lemma \ref{lem1} we obtain

\begin{align*}
I(m)=&\int_{0}^{1}\frac{1}{u^2+1}\bigg[\frac{m+1}{2}\log 2+\sum_{k=0}^{\frac{m-3}{2}}\log(1+\varphi_k)-m\log(1+u)\\
&+\sum_{k=0}^{\frac{m-3}{2}}\log\left(u^2+\frac{1-\varphi_k}{1+\varphi_k}\right)\bigg]{\rm d}u
\end{align*}
or

\begin{align*}
I(m)=&\frac{(m+1)\pi}{8}\log 2+\frac{\pi}{4}\sum_{k=0}^{\frac{m-3}{2}}\log(1+\varphi_k)-m\int_{0}^{1}\frac{\log(1+u)}{1+u^2}{\rm d}u\\
&+\sum_{k=0}^{\frac{m-3}{2}}\int_{0}^{1}\frac{\log\big(u^2+q\big)}{u^2+1}{\rm d}u,
\end{align*}
where $q=\frac{1-\varphi_k}{1+\varphi_k}$.
Since
$$
\int_0^1 \frac{\log (u+1)}{u^2+1} \, du=\frac{ \pi}{8}\log 2
$$
we get
\begin{align*}
I(m)=\frac{\pi}{8}\log 2+\frac{\pi}{4}\sum_{k=0}^{\frac{m-3}{2}}\log(1+\varphi_k)+\sum_{k=0}^{\frac{m-3}{2}}\int_{0}^{1}\frac{\log\big(u^2+q\big)}{u^2+1}{\rm d}u,
\end{align*}
Applying Lemma \ref{lem3} with $q=\frac{1-\varphi_k}{1+\varphi_k}$, we get
\begin{align*}
I(m)=&\frac{\pi}{8}\log2-\frac{(m-1)\operatorname{G}}{2}+\frac{\pi}{4}\sum_{k=0}^{\frac{m-3}{2}}\log(1+\varphi_k)\\
&+\frac{\pi}{2}\sum_{k=0}^{\frac{m-3}{2}}\log\bigg(1+\sqrt{\frac{1-\varphi_k}{1+\varphi_k}}\thinspace\bigg)+\sum_{k=0}^{\frac{m-3}{2}}\operatorname{Ti}_2\bigg(\frac{\sqrt{1-\varphi_k^2}-1}{\varphi_k}\bigg).
\end{align*}
Simplifying this expression we arrive at the desired result stated in Theorem \ref{thm2.2}.
\end{proof}


\begin{thm}\label{thm2.3}
Let the parameters $a$, $b$, $\alpha$, $\beta$, and the functions $\phi$ and $\varphi$ are such that the following integrals are convergent.
\begin{equation*}
\mathcal{F}(x;\alpha,\beta)=\alpha\int_{0}^{1}\phi^\prime(\alpha t)\phi(\beta t^x){\rm d}t
\end{equation*}
and
\begin{equation*}
\mathcal{G}(x;a,b)=a\int_{0}^{\infty}\varphi^\prime(at)\varphi(bt^x){\rm d}t.
\end{equation*}
Then we have
\begin{equation}\label{e.1.3}
\mathcal{F}(x;\alpha,\beta)+\mathcal{F}(1/x;\beta, \alpha)=\phi(\alpha)\phi(\beta)-\phi^2(0)
\end{equation}
and
\begin{equation}\label{e.1.4}
\mathcal{G}(x;a,b)+\mathcal{G}(1/x;b, a)=\lim\limits_{t\to\infty}\varphi(at)\varphi(bt)-\varphi^2(0).
\end{equation}
\end{thm}
\begin{proof}
By partial integration we have

\begin{align*}
\mathcal{F}(x;\alpha,\beta)=\phi(\alpha t)\phi(\beta t^x)\big|_{t=0}^{t=1}-\beta\int_{0}^{1}\phi(\alpha t)xt^{x-1}\phi^\prime(\beta t^x){\rm d}t.
\end{align*}
Making the change of variable $u=t^x$ it follows that
\begin{align*}
\mathcal{F}(x;\alpha,\beta)&=\phi(\alpha)\phi(\beta)-\phi^2(0)-\beta\int_{0}^{1}\phi^\prime(\beta u)\phi(\alpha u^{1/x}){\rm d}u\\
&=\phi(\alpha)\phi(\beta)-\phi^2(0)-\mathcal{F}(1/x;\beta, \alpha),
\end{align*}
which proves \eqref{e.1.3}. The proof of \eqref{e.1.4} can be achieved similarly.
\end{proof}

\section{Examples}
\begin{exam} \label{ex1}Letting $m=3$ in Theorem \ref{thm2.2} we obtain
\begin{equation}\label{e3.1}
\int_{0}^{1}\frac{\log\big(x^3+1\big)}{x^2+1}{\rm d}x=-\frac{\pi}{8}\log2-G+\frac{\pi}{2}\log\big(1+\sqrt{3}\big) -\operatorname{Ti}_2\big(2-\sqrt{3}\big).
\end{equation}
\end{exam}

\begin{rem}
Combining the formulas \eqref{e08.1} and \eqref{e3.1} we get
\begin{align}\label{e3.2}
\operatorname{Ti}_2\big(2-\sqrt{3}\big)&=\frac{2}{3}\operatorname{G}-\frac{\pi}{12}\log \big(2+\sqrt{3}\big).
\end{align}
The inverse tangent integral satisfies $\operatorname{Ti}_2(x)-\operatorname{Ti}_2\big(\frac{1}{x}\big)=\frac{\pi}{2}\log x$. Setting $x=2-\sqrt{3}$ here gives
\begin{align}\label{e3.3}
\operatorname{Ti}_2\big(2+\sqrt{3}\big)=\frac{2}{3}\operatorname{G}-\frac{\pi}{12} \log \left(362-209 \sqrt{3}\right).
\end{align}
Identity \eqref{e3.2} appears in \cite[Eq. (3.239), p. 216]{Valean1} and \eqref{e3.3} recovers the formula given in \cite[Eq. (13)]{Campbell}.
\end{rem}
\begin{exam}Letting $m=5$ in Theorem \ref{thm2.2} we get:
\begin{align*}
&\int_{0}^{1}\frac{\log\big(x^5+1\big)}{x^2+1}{\rm d}x=-\frac{3 \pi}{8}  \log 2-2\operatorname{G}+\operatorname{Ti}_2\left(1+\sqrt{5}-\sqrt{5+2 \sqrt{5}}\thinspace\right)\\
&+\operatorname{Ti}_2\left(1-\sqrt{5}+\sqrt{5-2 \sqrt{5}}\right)+\frac{\pi}{2}  \log \left(1+\sqrt{5}+\sqrt{10+2 \sqrt{5}}\thinspace\right).
\end{align*}
Here we used
\begin{equation*}
\cos \frac{\pi}{5}=\frac{1+\sqrt{5}}{4} \quad \mbox{and} \quad \cos \frac{3\pi}{5}=\frac{1-\sqrt{5}}{4}.
\end{equation*}
\end{exam}

\begin{exam} Setting $m=3$ in \eqref{e1} we get
\begin{align}\label{e3.4}
\int_{0}^{1}\frac{\log\big(x^3+1\big)}{x+1}{\rm d}x=\frac{1}{2}\operatorname{Li}_2\left(\frac{3}{4}\right)+\frac{3}{2}\log^22-\frac{\pi^2}{9}.
\end{align}
\end{exam}

\begin{exam}Setting $m=5$ in \eqref{e1} we get
\begin{align*}
&\int_{0}^{1}\frac{\log\big(x^5+1\big)}{x+1}{\rm d}x=\frac{5}{2}\log ^22-\frac{\pi^2}{5}+\frac{1}{2}\operatorname{Li}_2\left(\frac{5+\sqrt{5}}{8}\right)+\frac{1}{2}\operatorname{Li}_2\left(\frac{5-\sqrt{5}}{8}\right).
\end{align*}
\end{exam}

\begin{exam} Example \ref{ex1} together with \cite[Theorem 2.10]{Bat-2} (with m=3) lead to the following series evaluation:

\begin{align*}
\sum_{n=1}^{\infty}\frac{2^n\sum_{k=1}^{n}\frac{(3/2)^k}{k}}{(2n+1)\binom{2n}{n}}=\frac{\pi}{2}\log 2+\frac{10}{3}\operatorname{G}-\frac{2 \pi }{3}  \log \left(2+\sqrt{3}\right).
\end{align*}

\end{exam}

\begin{exam}Letting  $\alpha=-v$, $\beta=-u$ and  $\phi(t)=\log(1+t)$ in \eqref{e.1.3}, we have
\begin{equation*}
  \int_{0}^{1}\frac{\log(1-ut^x)}{v^{-1}-t}{\rm d}t+  \int_{0}^{1}\frac{\log(1-vt^{1/x})}{u^{-1}-t}{\rm d}t=-\log(1-u)\log(1-v).
\end{equation*}
\end{exam}
\begin{rem}
This gives a new, very short and elementary proof of \cite[Theorem 2.3 (1)]{Choie}.
\end{rem}
\begin{exam}For $\alpha=\beta=1$ and $\phi(t)=\log(1+t)$, we get from in \eqref{e.1.3}
\begin{equation*}
  \int_{0}^{1}\frac{\log(1+t^x)}{1+t}{\rm d}t+\int_{0}^{1}\frac{\log(1+t^{1/x})}{1+t}{\rm d}t=\log^22.
\end{equation*}
\end{exam}

\begin{exam} For $\alpha=\beta=1$ and $\phi(t)=\arctan(t)$, we get from in \eqref{e.1.3}
\begin{equation*}
  \int_{0}^{1}\frac{\arctan\big(t^x\big)}{1+t^2}{\rm d}t+  \int_{0}^{1}\frac{\arctan\big(t^{1/x}\big)}{1+t^2}{\rm d}t=\frac{\pi^2}{16}.
\end{equation*}
\end{exam}

\begin{exam} Setting  $\alpha=\beta=1$ and $\phi(t)=\arcsin(t)$ in \eqref{e.1.3}, we have
\begin{equation*}
  \int_{0}^{1}\frac{\arcsin\big(t^x\big)}{\sqrt{1-t^2}}{\rm d}t+  \int_{0}^{1}\frac{\arcsin\big(t^{1/x}\big)}{\sqrt{1-t^2}}{\rm d}t=\frac{\pi^2}{4}.
\end{equation*}
\end{exam}

\begin{exam} Setting in \eqref{e.1.3}, we have $\alpha=\beta=1$ and $\phi(t)=\mbox{arcsinh}(t)$  in \eqref{e.1.3}, we have
\begin{equation*}
  \int_{0}^{1}\frac{\mbox{arcsinh}\big(t^x\big)}{\sqrt{1+t^2}}{\rm d}t+ \int_{0}^{1}\frac{\mbox{arcsinh}\big(t^{1/x}\big)}{\sqrt{1+t^2}}{\rm d}t=-\log^2\big(1+\sqrt{2}\big).
\end{equation*}
\end{exam}

\begin{exam} For $\alpha=\beta=1$ and $\phi(t)=\operatorname{Li}_2(t)$, we get from \eqref{e.1.3}
\begin{equation*}
  \int_{0}^{1}\frac{\log(1-t)\operatorname{Li}_2\big(t^x\big)}{t}{\rm d}t+ \int_{0}^{1}\frac{\log(1-t)\operatorname{Li}_2\big(t^{1/x}\big)}{t}{\rm d}t=-\frac{\pi^4}{36}.
\end{equation*}
\end{exam}
\begin{exam} Setting $a=b=1$ and $\varphi(t)=\arctan(t)$ in \eqref{e.1.4}, we get
\begin{equation*}
\int_{0}^{\infty}\frac{\arctan(t^x)}{t^2+1}{\rm d}t+\int_{0}^{\infty}\frac{\arctan(t^{1/x})}{t^2+1}{\rm d}t=\frac{\pi^2}{4}.
\end{equation*}
\end{exam}
\vskip 3mm
\begin{flushleft}
\textbf{\large{Acknowledgements}}
\end{flushleft}
We would like to thank the editor and the reviewer for their very careful reading and numerous helpful suggestions.\\
\\
\textbf{Author contributions} Both authors have contributed equally in all aspects of the preparation  of this submission.\\
\\
\textbf{Funding} The authors have not disclosed any funding.\\
\\
\textbf{Data Availibility Statement} This manuscript has no associated data.\\
\\
\textbf{Declarations Conflicts of interest} The authors declare that they have no conflict of interest.
\bibliographystyle{amsplain}

\begin{thebibliography}{9}

\bibitem{Bat-2} Batır, N.: Logarithmic integrals with applications to BBP and Euler-type sums. Bull. Malays. Math.
Sci. Soc. 46 (96) (2023), 35pp.

\bibitem{Campbell} Campbell, J. M.: Special values of Legendre’s chi-function and the inverse tangent
integral. Irish Math. Soc. Bulletin. 89, (2022) 17–23.

\bibitem{Choie} Choie, Y. J. and Kumar, R.: Arithmetic properties of the Herglotz-Zagier-Novikov function. Adv. Math. 433 (2023) 109313.
\bibitem{Grad} Gradshteyn, I. S. and Ryzhik, I. M.: Table of Integrals, Series and Products. 7th Edition, Elsevier, 2007.

\bibitem{Nantan} Dasireddy, N. S.: A solution to another interesting sum involving classical harmonic number and central binomial coefficient,
Math. Student. 93 (3-4) 2024 106-114.

\bibitem{Herglotz} Herglotz, G.: Über die Kroneckersche Grenzformel für reelle quadratische Körper I. Ber. Verh. Sächs.
Akad. Wiss. Leipz. 75 (1923) 3–14.

\bibitem{Muzaf} Muzaffar, H. and Williams, K. S.: A restricted Epstein zeta function and the evaluation of some definite integrals. Acta Arith.
104 (1) (2002) 23–66 .

\bibitem{Novikov} Novikov, A. P.: Kronecker’s limit formula in a real quadratic field. Math. USSR. Izv. 17 (1981) 147.

\bibitem{Zagier} Radchenko, D. and Zagier, D.: Arithmetic properties of the Herglotz function. J. Reine Angew. Math. 797 (2023) 229-253.

\bibitem{Ramanujan} Ramanujan, S.: On the integral $\int_{0}^{x}\frac{\tan^{-1}(t)}{t}{\rm d}t$. J. Indian Math. Soc. 7 (1915) 93–96.

\bibitem{Rogers} Rogers, M. D.: A study of inverse trigonometric integrals associated with three-variable Mahler measures,and some  related identities, J. Number Theory,  121 (2006) 265-304.

\bibitem{So-Bat} Sofo, A. and Bat{\i}r, N.: Parametrized families of polylog integrals. Constr. Math. Anal. 4 (4) (2021) 400–419.

\bibitem{Valean1} Vălean, C. I.: (Almost) Impossible Integrals, Sums and Series. Springer, 2019.

\bibitem{Valean2} Vălean, C. I.:  More (Almost) Impossible Integrals, Sums and Series. Springer, Cham, 1st Edition, 2023.

\bibitem{Zagier 2} Zagier, D.: A Kronecker limit formula for real quadratic fields. Math. Ann. 213 (2), 153–184 (1975).


\end{thebibliography}

\end{document}